\documentclass[11pt, reqno]{amsart}
\usepackage{geometry}
\usepackage{amsmath, amsthm, amssymb, bm, color, framed, graphicx, hyperref, mathrsfs}
\usepackage{cite}
\usepackage[all]{xy}
\usepackage{bbm}
\usepackage{cancel}
\usepackage{esint}
\definecolor{shadecolor}{RGB}{241, 241, 255}

\newtheorem{theorem}{Theorem}

\newtheorem{lemma}{Lemma}

\def\al{\aligned}
\def\eal{\endaligned}
\def\be{\begin{equation}}
	\def\ee{\end{equation}}
\def\befr{\begin{frame}}
	\def\enfr{\end{frame}}

\def\a{\alpha}
\def\b{\beta}
\def\e{\epsilon}
\def\n{\nabla}

\def\R{{\bf R}}
\def\M{{\bf M}}

\def\p{\partial}
\def\de{\delta}
\def\De{\Delta}

\usepackage{amssymb,amsthm,amsmath,amsfonts}
\usepackage{xcolor,paralist,hyperref,fancyhdr,etoolbox}

\usepackage{mathtools}
\usepackage{esint}
\usepackage{colonequals}
\newcommand\dif{\mathop{}\!\mathrm{d}}

\def\Ric{\mathrm{Ric}}
\usepackage{microtype}

\makeatletter
\let\originalheadings\ps@headings
\def\ps@headings{%
  \originalheadings
  \def\@oddhead{%
    \normalfont\scriptsize
    \hbox to\textwidth{%
      \hfil
      \parbox[c]{.83\textwidth}{%
        \centering\rightmark
      }%
      \hfil\llap{\thepage}%
    }%
  }%
}
\makeatother
\title[Finiteness of fundamental groups for manifolds\protect\\
     with sign-changing Ricci curvature]
      {Finiteness of fundamental groups for manifolds with sign-changing Ricci curvature}
\begin{document}
\author{Lang Qin}
\address{School of Mathematics, Fudan University, Shanghai, China }
\email{lqin22@m.fudan.edu.cn}
\date{}	
\maketitle
\begin{abstract}
	 Let $\M$ be a complete open $\kappa$ noncollapsed Riemannian manifold and have Ricci curvature bounded from below. If the negative part of its Ricci curvature goes to zero sufficiently rapidly at infinity in terms of the volume noncollapse constant and dimension $n$, then the fundamental group of $\M$ is finite. When the Ricci  curvature is nonnegative, the result agrees with the classical conclusion of \cite{Li86}.
\end{abstract}
 \section{Introduction}
  Let $(\M^n,g)$ be an open Riemannian $n$-dimensional manifold. Throughout the paper, $V(x,R)$ denotes the volume of the geodesic ball of radius $R$ with center at $x$ in $\M$, $\mu=\mu_g$ denotes the Riemannian volume measure.   
  We  use the Laplace–Beltrami operator $\De=\mathrm{div} \n $. The heat equation is
  \be \p_t u =\De u. \ee
  When Ricci curvature is bounded below, it is well-known that there exists a unique  minimal heat kernel   $H(x,y,t)$(see \cite{D83} and \cite{Li}).
 
 The connection between Ricci curvature and the fundamental group is
 important to the study of open manifolds. Milnor showed that if $\M$ has nonnegative Ricci curvature, then any finitely generated subgroup of the fundamental group  has polynomial  growth, and conjectured that the fundamental group is finitely  generated \cite{M68}. Milnor's conjecture has been disproved by Bru\'e, Naber and Semola in dimensions $n\ge 6$ \cite{BNS25a,BNS25b}, and in dimensions $n=4,5$ recently by Wu and Yan \cite{WY26}. See \cite{CV35,SY82,Liu13,Pan20b}  for positive results in lower dimensions. In fact, these works establish stronger
 rigidity results. 
 
 Additional restrictions are added to study the fundamental group under nonnegative Ricci curvature.
 Li's work \cite{Li86} is one of the earliest contributions on the subject.
 \begin{theorem}[P. Li]
 	Let $(\M^n,g)$ be a complete noncompact 
 	Riemannian manifold with nonnegative Ricci curvature and maximal volume growth. Then
 	\be V(q,\sqrt t) H(q,q,t) \to (4\pi)^{-n/2}\omega_n, \ee for every fixed $q\in\M$, as $t$ goes to infinity. Here the constant $\omega_n$ is the volume of the unit ball in $\R^n$.
 	Moreover,  $\M$ has a finite fundamental group.
 \end{theorem}
 Anderson  studied manifolds with polynomial volume growth of order $\ge n-k$ 
and the first Betti number $\ge k$ to show they have finitely generated fundamental group \cite{A90}. Sormani proved finite
generation under a small linear diameter growth condition \cite{S00}.
Pan and Rong obtained finite generation results with an
additional condition on the deck action \cite{PR18}  and conjectured that adding Euclidean volume growth on the universal cover implies finite generation.  Huang and Huang established the $4$-dimensional case and Hongzhi Huang established the $5$-dimensional  case \cite{HH25,H26}. For more results on fundamental groups of manifolds with nonnegative Ricci curvature, we refer to \cite{G78,FY92,Wei88,KW11,W00,Pan20a},etc.

There is also a long history of results describing  how curvature decay and volume growth conditions restrict the topology. Abresch and Gromoll proved that manifolds satisfying \eqref{07} and \eqref{08} have finite topological type if the
sectional curvatures are uniformly bounded and the diameter growth has order of $o(r^{1/n})$ \cite{AG90}. Sha and Shen obtained finite topological type under
global nonnegative Ricci curvature, maximal volume growth, and a
quadratic lower sectional curvature bound \cite{SS97}. Lott and Shen
proved homotopy equivalence to a finite CW complex under a quadratic
lower sectional curvature bound, subquadratic volume growth, and
uniform noncollapsing of unit balls \cite{LS00}. We refer to \cite{M08,Z15} and the references within.

In this paper, we consider $\kappa$ noncollapsed manifolds whose Ricci curvature has faster than quadratic decay, and prove the following theorem:
\begin{theorem}
	Let $(\M^n,g)$ be a complete noncompact Riemannian manifold of dimension $n\ge 2$.  We fix a base point $p\in\M$, and denote by $$r=r(x)=d(x,p)$$ the distance function to $p$. Suppose that $\M$ is $\kappa$ noncollapsed in all scales:
	\be V(x,R) \ge \kappa R^{n},\label{03} \ee
	for some constant $\kappa>0$, for all $x\in\M$ and $R>0.$  Assume also that the Ricci
	curvature of $\M$ satisfies
	\be \Ric \ge 0\ \  \mbox{\ on\ }B(p,1).\label{02} \ee
	There exists some $\epsilon_0=\e_0(n,\kappa)>0$ such that, if \be \Ric \ge - \e^2 r^{-(2+\e)},\label{01} \ee 
	for some $0\le\e\le \e_0 $, then the fundamental group $\pi_1 (\M)$ is finite.
	\label{main}
\end{theorem} 

When $\e=0$, Ricci curvature is nonnegative, 
$\kappa$ noncollapsedness \eqref{03} is, by Bishop's comparison theorem,  automatically equivalent to maximal volume growth. Theorem \ref{main}. recovers Li's Theorem.  Volume growth
at one fixed center does not  provide the uniform bound
required in our heat kernel argument. It
is enough to assume \eqref{03} for every center and $R\ge 1$: Let $K=\e^2/(n-1)$, and let $v_{-K}(s)$ be the volume of the ball in the simply connected space form of curvature $-K$. The Bishop-Gromov inequality gives, for $0<R\le 1$, when $\e\le 1$,
$$V(x,R)\ge V(x,1) v_{-K}(R)/v_{-K}(1)\ge  \kappa \omega_n R^n/v_{-K}(1) \ge c_n\kappa R^n.$$  the coefficient is bounded below by a positive multiple
of $\kappa$ depending only on $n$,  the corresponding $\e_0^\prime=\min\{1,\e_0(n,c_n\kappa)\}.$

Condition \eqref{02} is a convenient choice to  control the compact region whose curvature is not controlled by \eqref{01}. Otherwise, consider the connected sum of $\R^n$ and a manifold with infinite fundamental group. Condition \eqref{01} can be replaced by \be \Ric \ge - (n-1)\lambda(r),\ \ \ r\ge 1, \label{07}\ee with  a suitable positive nonincreasing radial function $\lambda(r)$, for example, provided\be\lambda(1)+\int_1^\infty r\lambda(r)\dif r\le \e_0(n,\kappa). \label{08}\ee 

The dependence on $\kappa$ is important, for example, consider $\R\times S^1 $ and the warped product metric
\be g_\e=\dif s^2+f_\e^2 (s) \dif \theta^2
\ee with
\be f_\e(s)=\left\{
\al &1,\ \ \ &|s|\le 2,\\
&1+\frac{\e^2}{1000} (|s|-2) e^{-\frac{1}{|s|-2}},\ \ \ &|s|>2.
\eal \right.
\ee 
With $p=(0,0)$, it satisfies 
\be \al &\Ric=0 \ \ &\mbox{ on } B(p,1);\\
&\Ric\ge -\e^2r^{-2-\e}g_\e;  &\\
& V(q,R)\ge \frac{\e^2}{16000}R^2 \ \ \ &\mbox{ for all } q,R>0.
\eal 
\ee

Li used heat kernel method in \cite{Li86}, and Xu established a relation between the limiting behavior of the heat kernel at infinity and the behavior of the heat kernels on tangent cones at infinity with renormalized measure \cite{X14}. To apply Li's idea of comparing the heat kernels on the universal cover and the base manifold to estimate the order of fundamental group, we prove the following large time behavior of heat kernel:

\begin{theorem}
	Let $(\M^n,g)$ be a complete noncompact  Riemannian manifold. Suppose that $\M$ is  $\kappa$ noncollapsed in all scales, and  $\Ric \ge 0$ on $B(p,1)$.  There exists some $\epsilon_0=\e_0(n,\kappa)>0$ such that, if  \be \Ric \ge - \e^2 r^{-(2+\e)}\ee on $\M$ with $0\le \epsilon\le \epsilon_0$, then as $t$ goes to infinity,  the following holds
	for every fixed $q\in\M$
	\be V(q,\sqrt t) H(q,q,t) \to (4\pi)^{-n/2}\omega_n. \ee The constant $\omega_n$ is the volume of the unit ball in $\R^n$. Moreover, along any path $\gamma (t)$ in $\M\times\M\times (0,\infty)$ of the form
	\be \gamma(t)=(x(t),y(t),t) \ee with 
	\be d^2(q,x(t))+d^2(q,y(t)) =o(t)\ee as $t$ goes to infinity, we have
	\be V(q,\sqrt t) H(x(t),y(t),t) \to (4\pi)^{-n/2}\omega_n. \ee
	
	Equivalently, \be t^{n/2} H(x(t),y(t),t) \to (4\pi)^{-n/2}/AVR. \ee \label{hk}
	Here $AVR=\lim\limits_{R\to\infty} \frac{V(q,R)}{\omega_n R^n}$ is the asymptotic volume ratio of $\M$.
\end{theorem}

The analytic tools we use in this paper build on works of Petersen-Wei \cite{PW97}, who developed an integral Laplace
comparison theory. In \cite{ZZ18} Zhang and Zhu proved a Li-Yau gradient bound for collapsing manifolds under an integral curvature condition. Recently Carron, Mondello, Tewodrose, Impera, Rimoldi and Veronelli established more heat kernel
estimates for open manifolds whose Ricci curvature is controlled in the Kato class in \cite{C19}, \cite{CMT25} and \cite{IRV25}. 

\section*{AI Disclosure}\addcontentsline{toc}{section}{AI Disclosure}
This manuscript of proof was finished before the author started using AI tools. The full proof is due to the author. To illustrate the necessity of the dependence on $\kappa$, we add the example constructed by OpenAI in the previous section. The rest of the article does not contain any AI written text.

The author retains full responsibility for the correctness of the final work.

 \section{Preliminaries}
We first use a standard volume comparison argument (see, for example, \cite{Li} Chapter 1,2 and 4) to prove  that the asymptotic volume ratio exists.
\begin{lemma}
	Let $(\M^n,g)$ be a complete noncompact $\kappa$ noncollapsed Riemannian manifold with $\Ric \ge - \e^2 r^{-(2+\e)}$ for some $\e>0$, and $\Ric\ge 0$ on $B(p,1)$. Then for every fixed $q\in\M$,
	\be AVR=\lim\limits_{R\to\infty} \frac{V(q,R)}{\omega_n R^n}\ee exists. \label{L1}
\end{lemma}
\begin{proof}
  In terms of polar normal coordinates at $p$, we write the Riemannian volume measure using volume density $J$:
  \be \dif\mu_g=J(\theta,r)\dif r  \dif \theta, \ee before the cut time, with $\theta$ in the unit tangent sphere, where $\dif \theta $ is the area element of the unit $(n-1)-$sphere.        Moreover \be \frac{J(\theta,r)}{r^{n-1}} \to 1 \ee as $r\downarrow 0$.  We denote the cut time by $c(\theta)$ and                                
  write \be \hat J(\theta,r)=J(\theta,r)^{1/(n-1)}. \ee If $x=(\theta,r)$ is away from $p$ and the cut locus of $p$, the Bochner formula and $\De r=\p_r(\log J)$ give
  
  \be \al \frac{\p^2\hat J}{\p r^2} (\theta,r)=&-\frac{|\n^2 r|^2}{n-1} \hat J(\theta,r)- \frac{\Ric(\frac{\p}{\p r},\frac{\p}{\p r})}{n-1}\hat J(\theta,r)+\frac{(\De r)^2}{(n-1)^2} \hat J(\theta,r) \\ \le& -\frac{\Ric(\frac{\p}{\p r},\frac{\p}{\p r})}{n-1}\hat J(\theta,r)\le a(r) \hat J(\theta,r),\eal\label{11} \ee
for $0<r<c(\theta)$ before the cut time and
  \be a (r)=\frac{\max\{1,r\}^{-2-\e} \e^2}{n-1}. \ee
  Let $\phi:[0,\infty) \to\R $ be the solution to
  \be \phi^{\prime\prime} (s)=a(s)\phi(s), \ \ \phi(0)=0,\ \ \phi^\prime(0)=1.\label{13}\ee 
  
  From the equation for $\phi$ and its initial conditions, we have $\phi(s)\ge s$, $\phi^\prime (s)\ge 1$ and $\phi^{\prime\prime}(s)\ge 0$. Hence
  \be \phi(s)=\int_0^s\phi^\prime\le s\phi^\prime(s),\ \ \frac{\dif}{\dif s}(\log \phi^\prime)(s)=a(s)\frac{\phi(s)}{\phi^\prime(s)}\le a(s) s. \ee
  Integrating this, we obtain
  \be 1\le \phi^\prime(s)\le e^{\int_0^\infty s a(s)\dif s} , \ \ \ \ s\le \phi(s) \le  e^{\int_0^\infty s a(s)\dif s}s. \label{16}\ee
  Thus $\phi^\prime(\infty)$ exists and \be \lim\limits_{R\to\infty} \left(\frac{\phi(Rs)}{R}\right)^{n-1}=\lim\limits_{R\to\infty}\left(\frac{\int_0^{Rs}\phi^\prime(t)\dif t }{R}\right)^{n-1}= (s\phi^\prime(\infty))^{n-1}. \ee 
  
  Define the model volume by \be V_{\phi}(R)=n\omega_n\int_0^R \phi(s)^{n-1}\dif s. \ee
  
  The dominated convergence theorem implies that
  \be \lim\limits_{R\to\infty} \frac{V_{\phi}(R)}{R^n}=\lim\limits_{R\to\infty}n\omega_n\int_0^1 \left(\frac{\phi(Rs)}{R}\right)^{n-1}\dif s =\omega_n \phi^\prime(\infty)^{n-1}. \label{modelvolume}\ee 
  
   Comparing \eqref{11} and \eqref{13}, we have
  \be \frac{\p}{\p r} \left(\frac{\p \hat J}{\p r} \phi-\hat J \frac{\p \phi}{\p r}\right)\le 0.\ee 
  The initial asymptotic $\hat J(r)=r+O(r^3), \phi(r)=r+O(r^3)$
  imply $\frac{\p \hat J}{\p r} \phi-\hat J \frac{\p \phi}{\p r}\to 0$ as $r$ decreases to zero. 
  Therefore \be \frac{\p (\hat J/\phi)}{\p r}\le 0, \ee and \be \frac{J(\theta,r)}{\phi(r)^{n-1}}\label{21} \ee is nonincreasing in $(0,c(\theta))$. Extend $J$ by zero for $r\ge c(\theta)$, and the quotient \eqref{21} remains nonincreasing.
  Consequently \be A(r):=\frac{1}{n\omega_n \phi(r)^{n-1}}\int_{S_p\M} J(\theta,r)\dif\theta,\ \ A(0):=1 \ee is nonincreasing and nonnegative, with $A(r)\le A(0)=\lim\limits_{r\to 0}\frac{J/r^{n-1}}{(\phi/r)^{n-1}}=1.$ It follows that
  \be \frac{V(p,R)}{V_\phi(R)}=\frac{\int_0^R A(r)\phi(r)^{n-1}\dif r}{\int_0^R \phi(r)^{n-1}\dif r} \ee is a nonincreasing function of $R$. Together with \eqref{03} and \eqref{16}, we deduce that
  \be\frac{\kappa}{\omega_n e^{(n-1)\int_0^\infty s a(s)\dif s}} \le  \frac{V(p,R)}{V_\phi(R)} \le 1. \label{vc}\ee Thus the limit of $ \frac{V(p,R)}{V_\phi(R)}$ exists as $R\to\infty$. Combining this with \eqref{modelvolume}, we obtain
  \be \lim\limits_{R\to\infty} \frac{V(p,R)}{R^n}= \lim\limits_{R\to\infty}\frac{V(p,R)}{V_\phi(R)} \lim\limits_{R\to\infty} \frac{V_\phi(R)}{R^n} \ee exists.
  
  For $R>d(p,q)$ large enough, we have
  \be V(p,R-d(p,q))\le V(q,R) \le V(p,R+d(p,q)).\ee
  Hence 	\be  AVR= \lim\limits_{R\to\infty} \frac{V(q,R)}{ \omega_n R^n}\ee
exists. The limit is independent of $q$ and positive.
  
\end{proof}
Theorem 2.4. of \cite{PSR} gives the following Laplacian comparison:
\begin{lemma}
		Under assumptions \eqref{02}-\eqref{01}, let $\phi$ be the function defined in the
		proof of  Lemma \ref{L1}. Then
		\be \De r\le (n-1)\frac{\phi^\prime(r)}{\phi(r)} \ee 
		holds pointwise on $\M\setminus (Cut(p)\cup\{p\})$, and weakly on all of $\M$. 
		Hence \be \De r \le\frac{n-1+o(1)}{r} \ee as $r$ goes to infinity.
		\label{L2}
\end{lemma}

 \section{Large Time Behavior of the Heat Kernel}
 Our goal in this section is to prove Theorem \ref{hk}. If $\e=0$, the result follows from Li's theorem. We henceforth assume $\e>0$.
 
  Define
  \be \Ric_-(x):=\max\{0,-\inf\limits_{v\in T_x\M,g_x(v,v)=1}\Ric_x(v,v)\}, \ee  
  and
  \be k_T(\M,g):=\sup\limits_{x\in\M} \int_0^T \int_{\M} H(x,y,s)\Ric_-(y)\dif \mu_y\dif s,  \ee 
  where $H(x,y,s)$ is the heat kernel of $\M$. 
  
  We also write $k_\infty(M,g)=\sup\limits_{T>0} k_T(M,g)$ and the forward heat semigroup $P_t=e^{t \Delta}$.
  
  The following result is from Proposition 3.16, and Proposition 3.3 in \cite{C19} localized by Proposition 4.2 and Remark 4.3 in \cite{CMT25}.
  \begin{theorem}[G. Carron, I. Mondello and D. Tewodrose]
  	Let $(\M^n,g)$ be a complete Riemannian manifold. Suppose there exists some $T>0$ satisfying 
  	\be k_T(\M,g)\le \frac{1}{16n}.\ee
  	Let $u: \M\times[0,T] \to \R $ be a positive solution of the heat equation 
  	\be u_t=\De u. \ee
  	Then, on $\M\times  (0,T]$, we have 
  	\be e^{-8\sqrt {nk_T}}\frac{|\n u|^2}{u^2} -\frac{1}{u} \frac{\p u}{\p t} \le \frac{n}{2t} e^{8\sqrt {nk_T}}, \label{ly} \ee 
  	and \be H(q,q,t)\le \frac{C(n)}{V(q,\sqrt t)}, \label{lyh}\ee
  	for some constant $C(n)=2^{ne^2/2}e^{e^2/2}$.
  	\label{C}
  \end{theorem}
  
  Under a smallness condition on the Kato constant, Impera, Rimoldi and Veronelli proved the following result, Theorem 1.1 in \cite{IRV25}:
  \begin{theorem}[D. Impera, M. Rimoldi and G. Veronelli]
  	Let $(\M^n,g)$ be a complete nonparabolic Riemannian manifold of dimension $n\ge 3$. Suppose that \be \limsup\limits_{R\to\infty} \frac{V(y,R)}{R^n}\ge \kappa>0,\ee
  	for some $y\in\M$. Assume that\be k_\infty(\M,g)< \frac{1}{n-2}.\ee Then there exists $C(k_\infty,n,\kappa)
  	= n [1-(n-2)k_\infty]^{4(n-1)/(n(n-2))} \kappa^{1/n}$, \be 
  	\left(\int_{\M} |h|(x)^{\frac{n}{n-1}} \dif\mu_x\right)^{\frac{n-1}{n}} \le  C(k_\infty,n,\kappa) ^{-1}  \int_{\M} |\n h| \dif\mu_x, \label{Sob}
  	\ee for every $h\in \mathcal C_c^\infty(\M).$
  	 \label{I}
  	  \end{theorem}
  In fact, we choose a smooth compact domain $\Omega$ containing the support of $h$ and apply  the original theorem to $\sqrt{h^2+\eta^2}$ on $\Omega$ and let $\eta$ go to zero. 
  
  \begin{proof}[Proof of Theorem \ref{hk}]
  	 Case I. $n\ge 3$.
  	  
  	Under our assumptions \eqref{01} and \eqref{02}, we have \be 0\le \Ric_-(x)\le \e^2 1_{\{r(x)\ge 1\}} r(x)^{-2-\e}. \ee
  	 We need to bound $k_T$ uniformly. Observe that $k_0=0$ and $k_T$ is continuous in $T$:
  	 \be 0\le |k_{T_1}-k_{T_2}|\le \e^2|T_1-T_2|.\ee Thus we can fix some $T>0$ with $k_T<\frac{1}{16n}$. Theorem \ref{C} and \eqref{03} give the heat kernel bound \eqref{lyh} on $ \M \times (0,T] $. Together with the semigroup Cauchy–Schwarz inequality, we obtain:
  	\be H(x,y,s)\le \left(H(x,x,s)H(y,y,s)\right)^{1/2}\le  \frac{C(n,\kappa)}{s^{n/2}}, \label{h2}\ee 
  	for $0<s\le T/2$, with $C(n,\kappa)=C(n)/\kappa$.
  	
  	Split the integral as follows when $s\ge 1$:
  	\be 
  	\al 
  	&\int_{\M} H(x,y,s)\Ric_-(y) \dif \mu_y\\
  	=&\int_{\{1\le r(y)\le \sqrt s\}} H(x,y,s)\Ric_-(y) \dif \mu_y +\int_{\{ r(y)> \sqrt s\}} H(x,y,s)\Ric_-(y) \dif \mu_y \\
  	\le& \frac{C(n,\kappa)}{s^{n/2}}C_1(n)\e^2\int_1^{\sqrt s} r^{n-3-\e} \dif r +\e^2 s^{-1-\e/2} \int_{\{ r(y)> \sqrt s\}} H(x,y,s)\dif \mu_y\\
  	\le& C(n,\kappa)C_1(n)\e^2s^{-1-\e/2}+\e^2 s^{-1-\e/2}\\
  	\le& C_1(n,\kappa)\e^2s^{-1-\e/2},
  	\eal
  	\ee
  	where for the first part we use \eqref{h2}, polar coordinates and volume comparison \eqref{16},\eqref{21} and \eqref{vc}, for the second part we use \eqref{01} and the property of the heat kernel, and in the last inequality we choose $\e<\min\{1,\frac{n-2}{2}\}$, $C_1(n)=e^{3/2}n\omega_n$ and $C_1(n,\kappa)=C(n,\kappa)C_1(n)+1$.
  	
  	We now integrate in time, splitting the time interval into two parts if  $T/2\ge 1$:
  	\be
  	\al
  	&\int_0^{T/2} \int_{\M} H(x,y,s)\Ric_-(y) \dif \mu_y \dif s\\
  	=&\int_0^{1} \int_{\M} H(x,y,s)\Ric_-(y) \dif \mu_y \dif s +\int_1^{T/2} \int_{\M} H(x,y,s)\Ric_-(y) \dif \mu_y \dif s\\
  	\le&\e^2\int_0^{1} \int_{\M} H(x,y,s) \dif \mu_y \dif s +\int_1^{\infty}C_1(n,\kappa)\e^2s^{-1-\e/2} \dif s\\
  	\le&\e^2+2C_1(n,\kappa)\e,
  	\eal 
  	\ee using $\Ric_-\le\e^2$ directly if $T/2\le 1$,
  	which tells us
  	\be k_{T/2}(\M^n,g)\le C_2(n,\kappa)\e.\ee
  	
  	The semigroup property of the heat kernel implies
  	\be 
  	\al
  	\int_0^T P_s \Ric_-\dif s=&\int_0^{T/2}P_s \Ric_-\dif s+\int_{T/2}^T P_s \Ric_-\dif s\\
  	=&\int_0^{T/2}P_s \Ric_-\dif s+P_{T/2}\int_0^{T/2} P_s \Ric_-\dif s\\
  	\le&2k_{T/2}.
  	\eal
  	\ee 
  	Taking the supremum over $\M$, we have
  	\be k_T(\M^n,g)\le 2k_{T/2}(\M^n,g). \ee
  	
  	Suppose that $T^*$ is the first time that $k_T$ reaches $\frac{1}{16n}$.  The preceding
  	calculation shows 
  	\be k_{T^*}\le 2k_{T^*/2}\le 2C_2(n,\kappa)\e. \ee
  	
  	Choosing $\e\le \frac{1}{64nC_2(n,\kappa)}$, we obtain $\frac{1}{16n}=k_{T^*}\le \frac{1}{32n}$, a contradiction. Thus no such finite $T^*$ exists.
  	
  	The Kato bound now holds for every $T$, so Theorem \ref{C}  applies globally. The same argument shows that
  	\be H(x,y,t)\le C(n,\kappa) t^{-n/2},\label{53} \ee
  	and
  	\be \sup\limits_{x\in\M} \int_0^\infty \int_{\M} H(x,y,s)\Ric_-(y)\dif \mu_y\dif s \le C_2(n,\kappa)\e .
  	\label{54}  \ee
  	
  	Choose $\epsilon_0=\min\{1,\frac{n-2}{2},\frac{1}{64nC_2(n,\kappa)}\}$.
  	
  The estimate \eqref{53} implies that 
  	\be \int_1^\infty H(x,x,t) \dif t\le C \int_1^\infty t^{-n/2} \dif t<\infty,\ \ \    n \ge 3.
  	 \ee 
  	 hence $\M$ is nonparabolic.
  	 
  	 We  apply Theorem \ref{I} to derive a global Sobolev inequality \eqref{Sob}:
  	\be 
  	\left(\int_{\M} |h|(x)^{\frac{n}{n-1}} \dif\mu_x\right)^{\frac{n-1}{n}} \le  C(k_\infty,n,\kappa) ^{-1}  \int_{\M} |\n h| \dif\mu_x .
  	\ee 
  	
  	The Sobolev inequality shows that for $u\in\mathcal C_c^\infty(\M)$
  	\be \left(\int_{\M} |u|^{2n/(n-2)} \dif\mu_x \right)^{(n-2)/n} \le C \int_{\M} |\n u|^2 \dif\mu_x, \label{56}\ee 
  	with $C=\left(\frac{2(n-1)}{(n-2)C(k_\infty,n,\kappa)}\right)^2$.

  	To determine the large time behavior of the heat kernel,   	we define 
  	\be f_R(x,t)= \frac{1}{AVR(4\pi t)^{\frac{n}{2}}R^n} \exp \left( - \frac{r(x)^2}{4R^2 t} \right). \ee 
  	
  	 Integration by parts and Lemma \ref{L1} show that
  	\be \al \int_\M f_R(x,t)\dif\mu_x=&\frac{1}{AVR(4\pi t)^{\frac{n}{2}}R^n} \int_0^{\infty } \exp \left( - \frac{s^2}{4R^2 t} \right)\dif V(p,s)\\ 
  	=&\frac{1}{AVR(4\pi t)^{\frac{n}{2}}R^n} \int_0^{\infty } \exp \left( - \frac{s^2}{4R^2 t} \right)V(p,s)\frac{s}{2R^2t}\dif s\\
  	=&\frac{1}{AVR(4\pi t)^{\frac{n}{2}}} \int_0^{\infty } \exp \left( - \frac{z^2}{4 t} \right)\frac{V(p,Rz)}{(Rz)^n}\frac{z^{n+1}}{2t}\dif z\to 1 \eal \label{58}\ee
  	for fixed $t>0$ and $R\to\infty$. In fact, the convergence is compact uniform for $t$,
  	\be \sup\limits_{\de\le t\le T} \left|\int_\M f_R(x,t)\dif\mu_x-1\right|\to 0,\ \ \ \sup\limits_{\de\le t\le T} \left|\p_t\int_\M f_R(x,t)\dif\mu_x\right|\to 0. \ee
  	
  	Away from $p$ and the cut locus, using Lemma \ref{L2}, we derive that
  	\be 
  	(\p_t-R^2\De)f_R =\frac{f_R}{4t}\left(\De (r(x)^2) -2n\right)\le \frac{f_R}{4t} e(r),
  	\label{59}
  	\ee 
  	where\be e(r)=2(n-1) [r \phi^\prime(r)/\phi(r) - 1] \ee
  	 is nonnegative, bounded, and tends to zero.
  	From now on, we regard $\De f_R$ as as a Radon measure (see \cite{EG}), and define
  	\be \Lambda(R,t):=(\p_t-R^2\De)f_R=(\p_t f_R)\mu_g-R^2\De f_R. \ee   \eqref{59} gives
  	\be \De (r^2) \le (2n+e(r)) \mu. \ee 
  	Consequently
  	\be
  	\sigma:=(2n+e(r)) \mu-\De (r^2)
  	\ee
  	is a nonnegative Radon measure. Since $r^2$ is smooth near $p$, $\sigma$ has no atom at $p$.
  	Moreover, 
  	\be \Lambda(R,t)= \frac{f_R}{4t} e(r) \mu-\frac{f_R}{4t} \sigma. \ee
  	In particular, \be [\Lambda(R,t)]_+\le \frac{f_R}{4t} e(r)\mu. \ee The singular cut locus
  	contribution belongs to the negative part of $\Lambda(R,t)$.
  	
  	Consider a Lipschitz cutoff function
  	\be \chi_L=\left\{\al &0 \mbox{ outside } B(p,2L)\\&1 \mbox{ on } B(p,L)\eal \right.
  	\ee   with $|\n\chi_L|\le\frac{C}{L}.$   For fixed $R,t$,
  	\be
  	R^2 \left|\int_\M \langle\n\chi_L, \n f_R\rangle\dif\mu\right|\le \frac{C}{t} \int_{L\le r(x)\le 2L} f_R \dif\mu.
  	\ee 	
  	
  	Polynomial volume growth and Gaussian decay give integrability.  We calculate that
  	\be 
  	\int_\M \chi_L (\p_t-R^2\De)f_R=\int_\M \chi_L \p_t f_R\dif\mu+R^2  \langle\n\chi_L, \n f_R\rangle \dif\mu -\int_{\p B(p,2L)\cap Cut(p)} \chi_L R^2\p_n f_R
  	\ee 
  	Letting $L\to\infty$,
  	\be \int_\M\dif \Lambda(R,t)=\int_{\M}(\p_t-R^2\De)f_R = \p_t\int_\M f_R(x,t)\dif\mu_x.\ee
  	
  	By monotone convergence applied to the positive and negative
  	parts, we deduce that $\Lambda(R,t)$ has finite total variation:

  	\be\al  \|\Lambda(R,t)\|_{\rm{TV}}
  	=&[\Lambda(R,t)]_+(\M)+[\Lambda(R,t)]_-(\M)=2[\Lambda(R,t)]_+(\M)-\Lambda(R,t)(\M)\\
  	\le& 2\int_{\M}\frac{f_R}{4t}e(r(x))\dif\mu_x-\p_t \int_{\M}f_R \dif\mu_x\\
  	\le& 2\int_{\{r(x)\le A\}}\frac{f_R}{4t}e(r)\dif\mu_x+2\sup\limits_{r(x)\ge A}e(r)\int_{\M}\frac{f_R}{4t}\dif\mu_x +\left|\p_t \int_{\M}f_R \dif\mu_x\right|.\eal \ee
  	Letting $R\to\infty$ and then $A\to\infty$, we have
  	\be
  \sup\limits_{\de\le t\le T} \|\Lambda(R,t)\|_{\rm{TV}}\to 0.
  	\ee
  	
  	 Polynomial volume	growth gives $f_R(\cdot,\tau)$ in $W^{1,2}(\M)$ and $\p_\tau f_R(\cdot,\tau)$  in $L^2(\M)$. Consider $\chi\in\mathcal C_c^\infty(\M)$, 
  	 \be \al
  	 &\frac{\dif}{\dif\tau} \langle f_R(\cdot,\tau),P_{R^2(t-\tau)}\chi \rangle_{L^2(\M)}\\
  	 =&
  	 \langle \p_\tau f_R(\cdot,\tau), P_{R^2(t-\tau)}\chi \rangle_{L^2(\M)}+
  	 R^2\int_{\M}  \langle  \n f_R(\cdot,\tau), \n (P_{R^2(t-\tau)}\chi)\rangle \dif\mu\\=&\int_{\M} P_{R^2(t-\tau)}\chi \dif \Lambda(R,\tau).\eal 
  	 \ee 
  		Integrating with respect to $\tau$ and repeating the cutoff argument, we derive a weak Duhamel formula for $f_R$ as an identity of finite signed measures 
  		\be f_R(\cdot,t)=P_{R^2(t-s)}f_R(\cdot,s)+\int_s^t P_{R^2(t-\tau)}\Lambda(R,\tau)\dif \tau. \ee
  		Here the heat semigroup acts on finite signed measures, which is measurable in $\tau$. The map $\tau\to \Lambda(R,\tau)$ is weakly measurable by its representation
  		in terms of the fixed measures $\mu$ and $\sigma$.
  		Hence
  		\be\|f_R(\cdot,1)-\int_\M H(\cdot,y,R^2(1-\de))f_R(y,\de)\dif\mu_y\|_{L^1}\le \int_\de^1 \|\Lambda(R,\tau)\|_{\rm{TV}}\dif\tau\to 0, \ee 
  		for every fixed $0<\de<1$ and $R\to\infty$.
  		
  		Using \eqref{58} and the $L^1$ contractivity of the heat semigroup, 
  		\be\left\|f_R(\cdot,1)-\int_\M H(\cdot,y,R^2(1-\de))\frac{f_R(y,\de)}{\int_\M f_R(x,\delta)\dif\mu_x}\dif\mu_y\right\|_{L^1}\to 0,\label{63}\ee 
  		for every fixed $0<\de<1$ and $R\to\infty$.
  		
  	According to the results of J. Nash and N. Varopoulos \cite{N58,V85}, the Sobolev inequality \eqref{56} is equivalent to a Gaussian upper bound on the heat kernel:
  	\be H(x,y,t)\le C t^{-\frac{n}{2}} \exp\left({-\frac{d^2(x,y)}{Ct}}\right).\label{64}\ee 
  	
  	We follow a standard heat kernel estimate argument from \eqref{53} and $\De P_t=P_{t/4}\De P_{t/2}P_{t/4}$,
  	\be |\p_t H|=|\De H|=|\langle H(x,\cdot,t/4), \De P_{t/2}H(y,\cdot,t/4)\rangle | \le Ct^{-n/2-1}.\label{8811}\ee 
  	Multiplying \eqref{ly} by $H^2$ and applying \eqref{53} and \eqref{8811}, we have
  	\be |\n_xH(x,y,t)|\le Ct^{-(n+1)/2}. \label{65}\ee
  	
  	Applying this estimate and \eqref{64} at $t=R^2(1-\de)$, for $L\ge 2$ and $0<\delta\le 1/2$ and $R>0$,  we calculate that
  	
  	\be\al &\left\|\int_\M H(\cdot,y,R^2(1-\de))\frac{f_R(y,\de)}{\int_\M f_R(x,\delta)\dif\mu_x}\dif\mu_y-H(p,\cdot,R^2(1-\de))\right\|_{L^1}
  	\\\le&\int_\M \|H(\cdot,y,R^2(1-\de))-H(p,\cdot,R^2(1-\de))\|_{L^1} \frac{f_R(y,\de)}{\int_\M f_R(x,\delta)\dif\mu_x}\dif\mu_y\\
  	\\\le&\int_{B(p,R)}\frac{f_R(y,\de)}{\int_\M f_R(x,\delta)\dif\mu_x}\left[\int_{B(p,LR)}\left|H\left(x,y,R^2(1-\delta)\right)-H\left(p,x,R^2(1-\delta)\right)\right|\dif\mu_x\right]\dif\mu_y\\+&\int_{B(p,R)}\frac{f_R(y,\de)}{\int_\M f_R(x,\delta)\dif\mu_x}\left[\int_{M\setminus B(p,LR)}\left(H\left(x,y,R^2(1-\delta)\right)+H\left(p,x,R^2(1-\delta)\right)\right)\dif\mu_x\right]\dif\mu_y\\+&2\int_{\M\setminus B(p,R)}\frac{f_R(y,\de)}{\int_\M f_R(x,\delta)\dif\mu_x}\dif\mu_y
  	\\\le& CL^n\int_{B(p,R)}\left(\frac{r(y)}{R}\right) \frac{f_R(y,\de)}{\int_\M f_R(x,\delta)\dif\mu_x}\dif\mu_y+Ce^{-c(L-1)^2}\int_{B(p,R)}\frac{f_R(y,\de)}{\int_\M f_R(x,\delta)\dif\mu_x}\dif\mu_y
  	\\+&2\int_{M\setminus B(p,R)}\frac{f_R(y,\de)}{\int_\M f_R(x,\delta)\dif\mu_x}\dif\mu_y
  	\\&\le CL^n\delta^{1/2}+Ce^{-c(L-1)^2}+Ce^{-c/\delta}.
  	\eal \label{66}\ee
  In the last inequality,	we used\be \int_{B(p,R)}\left(\frac{r(y)}{R}\right) f_R(y,\de)\dif\mu_y\le C\de^{1/2}. \ee

  Now fix $0<\eta\le 1$ and take $x\in B(p,\eta R)$. 
  \be R^n |H(p,x,R^2(1-\de))-H(p,p,R^2(1-\de))|\le R^nC(R^2(1-\de))^{-n/2} \frac{r(x)}{R\sqrt{1-\de}}\le C\eta.\label{8855}\ee
  	
  	On the other hand,
  	\be |R^{n} f_R(x,1) - (4\pi)^{-n/2}/AVR|= (4\pi)^{-n/2}/AVR|e^{-r(x)^2/4R^2}-1|\le C\eta^2.\ee
  	
  	Together, we have
  	
  	\be |R^{n} H(p,p,R^2(1-\de)) - (4\pi)^{-n/2}/AVR|\le C\eta + R^n|f_R(x,1)-H(p,x,R^2(1-\de))|.\ee
  	
  	We integrate over $B(p,\eta R)$ and divide by its volume:
  	\be\al  &|R^{n} H(p,p,R^2(1-\de)) - (4\pi)^{-n/2}/AVR|\\\le& C\eta+ \frac{ R^n\int_{B(p,\eta R) }|f_R(x,1)-H(p,x,R^2(1-\de))|\dif\mu_x}{V(p,\eta R)}\\
  	\le & C\eta+\frac{1}{\kappa\eta^n} \left\|f_R(\cdot,1)-H(p,\cdot,R^2(1-\de))\right\|_{L^1} .\eal \ee
  	
  	Taking the limit shows
  	\be\limsup\limits_{R\to\infty}|R^{n} H(p,p,R^2(1-\de)) - (4\pi)^{-n/2}/AVR|\le C\eta+\frac{CL^n\delta^{1/2}+Ce^{-c(L-1)^2}+Ce^{-c/\delta}}{\kappa\eta^n} .\ee
  	For fixed $\de$ in $(0,1/2)$, substitute $t=R^2(1-\de)$. 
  	Letting   $\de\to 0$ first, then $L\to \infty$, and finally $\eta\to 0$, we will have the desired limit.
  	By symmetry of the heat kernel and \eqref{65}, integration along minimizing geodesics in
  	each spatial variable suggests
  	
  	\be t^{n/2}|H(q,q,t)-H(p,p,t)| \le 2C \frac{d(p,q)}{\sqrt t}\to 0\ee
  	and
  	\be \limsup\limits_{t\to\infty} t^{n/2} |H(x(t),y(t),t) - H(q,q,t)|
  	\le\lim\limits_{t\to\infty} C \frac{d(q,x(t)) + d(q,y(t))}{\sqrt t}=0.
\label{9911}  	 \ee
  	
  	Case	II. For $n=2$, the curvature condition $\eqref{01}$ becomes:
  		$K_-\le \e^2r(x)^{-2-\e}$ where $K$ is the sectional curvature. 
  		The proof in Section 4, which uses only the geometric
  		assumptions, shows for sufficiently small $\e$ depending on $\kappa$ that the
  		connected orientation cover is diffeomorphic to $\R^2$. The covering has degree
  		at most two.  Since a finite proper quotient of a one-ended space  has one end, $\M$ has one end.
  		Direct calculation shows $\M$ belongs to class (a) of \cite{GS05} 	Section 7.5. By  Corollary 7.12 in \cite{GS05}, we have the global
  		parabolic Harnack inequality. Theorem 2.7 of \cite{GS05} gives
  		\be
  		H(x,y,t)\le C t^{-1}\exp\left(-\frac{d(x,y)^2}{Ct}\right).  		\ee
  		There are constants $C$ and $\a\in(0,1]$, such that
  		\be |H(x,y,t)-H(x',y,t)|\le C t^{-1} \left(\frac{d(x,x')}{\sqrt t} \right)^\a \ee
  		whenever $d(x,x')\le \sqrt t$, uniformly in all points and $t>0$.
  		The  weak Duhamel arguments above are dimension
  		independent.  Replace $C \eta$ in \eqref{8855} by $C\eta^\a$, and respectively $(\frac{r(y)}{R})^\alpha,\delta^{\alpha/2},(\frac{d(p,q)}{\sqrt t})^\alpha$ in\eqref{66}-\eqref{9911}. Everything else is the same as above.
  \end{proof}
  
  \section{Finiteness of Fundamental Group}
  Unlike the case in \cite{Li86}, we cannot lift the curvature assumption \eqref{01} directly to the universal cover $\widetilde{\M}$. One approach uses remote balls, see  \cite{GS05} for the definition. In this paper, we apply a consequence of Carron's Theorem 1.4 and Lemma 1.14 \cite{C19} to
  the universal cover $\widetilde{\M}$. 
  \begin{theorem}[G. Carron]
  	Let $(\M^n,g)$ be a complete Riemannian manifold of dimension $n\ge 3$. Suppose that
  	(a). $(\M^n, g)$ satisfies the Euclidean Sobolev inequality with
  	Sobolev constant $\a$:
  	\be \forall f\in \mathcal C_c^{\infty}(\M):\ \ \  \a \left(\int_{\M} |f|^{2n/(n-2)} \dif\mu_g \right)^{(n-2)/n} \le \int_{\M} |\n f|^2 \dif\mu_g. \ee
  	(b). There is a constant $\b>0$ such that
  	\be \sup\limits_{x\in\M}  \int_{\M} \int_0^{\infty} H(x,y,s)\dif s\Ric_-(y)\dif \mu_y<\frac{1}{(1+\b)(n-2)}. \ee 
  	Then there exists $c>0$ depending only on $n,\a$ and $\b$, such that for all $x\in\M$ and $R\ge 0$:
  	\be c^{-1} R^n\le V(x,R)\le cR^n.\label{CV} \ee 
  	\label{C2}
  \end{theorem}
  
  Now we are ready to prove the main theorem.

  \begin{proof}[Proof of Theorem \ref{main}]
  	{\it Case I. $n\ge 3$. }
  	We want to apply Theorem \ref{C2} to $\widetilde\M$. Let $\pi:\widetilde\M \to\M$ be the universal cover. $\tilde g=\pi^* g,\tilde\n=\n_{\tilde g},\tilde\mu=\mu_{\tilde g},\widetilde\Ric=\pi^*\Ric, \widetilde\Ric_-=\Ric_-\circ\pi$ are the lifted metric, connection, volume measure and Ricci tensor from the underlying $g$. First we check condition (a). Take any $f\in \mathcal C_c^\infty (\widetilde{\M})$, and define \be u(x):=\left(\sum\limits_{\tilde x\in \pi^{-1}(x)}|f(\tilde x)|^2\right)^{1/2}, \label{81}\ee
  	the sum is locally finite because $f$ has compact support.
  	
  	By the fact that $\pi$ is a local isometry and Cauchy-Schwarz, we have
  	\be |\n u(x)|^2\le \sum\limits_{\tilde x\in \pi^{-1}(x)}|\tilde\n f(\tilde x)|^2 \ee almost everywhere and since $2n/(n-2)>2$,
  	\be u(x)^{2n/(n-2)}\ge \sum\limits_{\tilde x\in \pi^{-1}(x)} |f(\tilde x)|^{2n/(n-2)} .\label{83}\ee

  	Since $u$ is compactly supported and locally Lipschitz, we can choose a sequence of $u_k\in\mathcal C_c^\infty(\M)$ converging to  $u$ in $W^{1,2}(\M)$.  After passing to a
  	subsequence, $u_k\to u$ almost everywhere.
  	 Applying
  	\eqref{56}  and Fatou's
  	lemma proves the same Sobolev inequality for $u$.  Integrating \eqref{81}-\eqref{83}  shows that
  	\be \a \left(\int_{\widetilde\M} |f|^{2n/(n-2)} \dif\mu_{\tilde g} \right)^{(n-2)/n} \le \int_{\widetilde\M} |\tilde\n f|^2 \dif\mu_{\tilde g}, \ee 
   where $\a$ depends on $n$ and $\kappa$.
  	
  	Condition (b) holds because of \eqref{54}, $\tilde P_t\pi^*\le \pi^*P_t$ and Tonelli's theorem:
  	\be\al  \int_{\widetilde\M} \int_0^{\infty} \widetilde H(\tilde x,\tilde y,s)\dif s\widetilde{\Ric}_{-}(\tilde y)\dif \tilde\mu_{\tilde y}=& \int_0^{\infty} \int_{\widetilde\M}  \widetilde H(\tilde x,\tilde y,s)\widetilde{\Ric}_{-}(\tilde y)\dif \tilde\mu_{\tilde y}\dif s
  	\\
  	\le& \int_0^{\infty} \int_{\M}  \sum\limits_{\tilde y\in \pi^{-1}(y)}\widetilde H(\tilde x,\tilde y,s){\Ric}_{-}( y)\dif \mu_{ y}\dif s\\
  	\le& \int_0^{\infty} \int_{\M}  H( \pi(\tilde x), y,s){\Ric}_{-}( y)\dif \mu_{ y}\dif s\\
  	\le&\sup\limits_{ x\in \M} \int_0^{\infty} \int_{\M}  H( x, y,s){\Ric}_{-}( y)\dif \mu_{ y}\dif s\\
  	\le& \frac{1}{16n}<\frac{1}{2(n-2)}. \eal\ee 
  	Finally, we take the supremum over $\tilde x$ in $\widetilde \M$.
  	
  	Fix $\b=1$ in Theorem \ref{C2}. Since  $\a$ depends only on $n$ and $\kappa$, the
  	constant $c$ in \eqref{CV} also depends only on $n$ and $\kappa$. Hence
  	\be \widetilde V(\tilde x, R) \le C_3(n,\kappa) R^n .\ee

  	Let $S$ be a finite nonempty subset of $\pi_1(\M)$, and let $p$ be the fixed
  	base point. Choose $\tilde p$ in $\pi^{-1}(p)$, and define
  	\be l(S):=\max\limits_{\gamma \in S} d(\tilde p, \gamma \tilde p). \ee
 \be  A:=\left\{\tilde\sigma (1)\Big| \al x\in \M\setminus(Cut(p)\cup\{p\}) ,\sigma \mbox{ is the unique minimizing geodesic from } p \mbox{ to } x,\\ \tilde\sigma \mbox{ is its unique lift starting from }\tilde p, \sigma(0)=p,\sigma(1)=x, \tilde\sigma(0)=\tilde p \eal\right \}.\ee
  Since the cut locus has
  	measure zero,  the restriction $\pi|_A$ is a measure preserving bijection onto $\M$ minus a measure zero, we have
  	\be \mathrm{Vol}_{\widetilde \M} (A\cap B(\tilde p, R))=V(p,R). \ee
  	
  	For every $\gamma\in S$, 
  	\be \gamma(A\cap B(\tilde p, R))\subset B(\tilde p, R+l(S)). \ee

  	Since the translates $\gamma(A)$, $\gamma$ in $S$, are pairwise disjoint,
  	\be |S|V(p,R)=\sum\limits_{\gamma\in S}\mathrm{Vol}_{\widetilde \M}(\gamma(A\cap B(\tilde p, R))) \le \widetilde V(\tilde p, R+l(S)). \ee
  	
  	Therefore, \be |S|\le \frac{\widetilde V(\tilde p, R+l(S))}{V(p,R)} \le C_3(n,\kappa)\frac{(R+l(S))^n}{V(p,R)}, \ee

  	Since $S$ is fixed and finite,  $l(S)$ is fixed. Letting
  	$R\to\infty$,  we have
  	\be |S|\le \frac{C_3(n,\kappa)}{\omega_n AVR}.\ee
  
  Consequently, $\pi_1(\M)$ is finite,  with order at most $\frac{C_3(n,\kappa)}{\omega_n AVR}.$

  	{\it Case II. $n=2$. } $\int_\M K_-\le 2\pi e^\e \e<\infty.$ Huber's theorem \cite{H57} gives finite topological type. Consider the orientation cover $\hat \M$ with degree at most two. Suppose $\pi_1(\hat \M)\not=0$, then  the fact that the covering is a local isometry and an inequality by \cite{F41,H64,ST89,ST93} (see (2.1)-(2.3) in \cite{BC14}) give\be L(b)-L(a)\le\int_a^b [2\pi\chi(B(\hat p, t))-\int_{B(\hat p,t)}\hat K]\dif t\le (b-a) \int_{\hat\M} \hat K_-,\ee where $L(t)$ is the length of boundary of $B(\hat p, t)$ and $\chi$ is Euler characteristic. Since $\frac{\dif}{\dif R}\hat V(\hat p , R)=L(R)$ almost everywhere \cite{BC14}. Choose $b=R,a=a_0$ large enough  for $B(\hat p,a_0)$ to have nonpositive Euler
  	characteristic and integrate again, \be \frac{\hat V(\hat p,R)}{R^2}\le C_0 R^{-2}+C_1 R^{-1}+ \frac{1}{2}\int_{\hat \M} \hat K_- . \ee 
  	Thus
  	\be \kappa\le \limsup\frac{\hat V(\hat p,R)}{R^2}\le \frac{1}{2}\int_{\hat \M} \hat K_-\le \int_{\M} K_-\le 2\pi e^\e \e<\kappa ,\ee a contradiction when choosing $\e$ small depending on $\kappa$. Therefore, $\pi_1(\hat \M)=0$ and $|\pi_1(\M)|\le 2$. In fact, $\hat\M$ must be  diffeomorphic to  $\R^2$. Since $1=\chi(\hat \M)=2\chi( \M)$ is impossible, $\pi_1(\M)$ must be trivial.
  \end{proof}
\section*{Acknowledgements}\addcontentsline{toc}{section}{Acknowledgements}
  	The author thanks Professor Bobo Hua  and  Professor Qi S. Zhang for  constant guidance and support.

\end{document}